\documentclass[twoside, 10pt]{article}
\usepackage{mathrsfs,amsfonts,amsmath}
\RequirePackage{ifthen,calc}
\usepackage{theorem}
\usepackage[dvips]{color}

\renewcommand{\baselinestretch}{1.0}

 \catcode`@=11
 \def\@evenhead{\hbox to\textwidth{\footnotesize\rm\thepage \hfill
  {\it Lu, P., Li, Y. and Yao Q.}}} 

 \def\@oddhead{\hbox to \textwidth{\footnotesize{\it
Record Numbers in Random walks} \hfill\thepage}}

 \renewcommand{\section}{\makeatletter
 \renewcommand{\@seccntformat}[1]{{\csname the##1\endcsname.}\hspace{0.45em}}
 \makeatother \@startsection
{section}
{1}
{0pt}
{0.25\baselineskip}
{0.25\baselineskip}
{\normalsize\bfseries\mathversion{bold}}}
\catcode`@=11

\newcommand\acks{\section*{Acknowledgements}}
\newtheorem{theorem}{Theorem}
\newtheorem{corollary}[theorem]{Corollary}
\newtheorem{lemma}[theorem]{Lemma}

\newtheorem{example}{Example}
\theoremheaderfont{\normalfont\bfseries}
\theorembodyfont{\slshape} {\theorembodyfont{\rmfamily}}
{\theorembodyfont{\rmfamily} }
{\theorembodyfont{\rmfamily}}
{\theorembodyfont{\rmfamily} }
{\theorembodyfont{\rmfamily} }
{\theorembodyfont{\rmfamily} } {\theorembodyfont{\rmfamily} } {\theorembodyfont{\rmfamily}
}
 \def\beqlb{\begin{eqnarray}}\def\eeqlb{\end{eqnarray}}
 \def\beqnn{\begin{eqnarray*}}\def\eeqnn{\end{eqnarray*}}
 \def\beqls{\begin{equation}}\def\eeqls{\end{equation}}
 \numberwithin{equation}{section}
\def\e{\mathrm{e}}

\def\P{\mathrm{P}}
\def\E{\mathrm{E}}
\def\e{\mathrm{e}}

\begin{document}

\title{\LARGE\bf Limit Properties of Record Numbers in Random walks}
\author{Penghui LU\renewcommand{\thefootnote}{\fnsymbol{footnote}}\footnotemark[1]\qquad Yuqiang LI\renewcommand{\thefootnote}{\fnsymbol{footnote}}\footnotemark[2]\qquad Qiang YAO\renewcommand{\thefootnote}{\fnsymbol{footnote}}\footnotemark[3]
\\ \small Key Laboratory of Advanced Theory and Application in Statistics and Data Science-MOE,
\\ \small School of Statistics, East China Normal University,
\\ \small Shanghai 200062, P. R. China.
}
\date{}

\maketitle
\renewcommand{\thefootnote}{\fnsymbol{footnote}}
\renewcommand{\thefootnote}{\fnsymbol{footnote}}
\renewcommand{\thefootnote}{\fnsymbol{footnote}}\footnotetext[3]{Corresponding author. Email address: qyao@sfs.ecnu.edu.cn}

\begin{abstract}

\noindent In this paper, we systematically summarize and enhance the understanding of weak convergence and functional limits of record numbers in discrete-time random walks under Spitzer's condition, and extend these findings to $\sigma$--record numbers using similar methods. Additionally, we identify a sufficient condition for the existence of functional limits for record numbers in continuous-time random walks. Finally, we derive corresponding results for large deviations, moderate deviations, and laws of the iterated logarithm pertaining to record numbers in discrete-time random walks.

\smallskip

\noindent {\bf Keywords}:\;{\small Discrete-time random walks; Continuous-time random walks; Record numbers; $\sigma$-Record Numbers; Weak convergence; Functional limits; Large and moderate deviations; Laws of the iterated logarithm.}

\smallskip

\noindent{\bf AMS 2000 Subject Classification:}\;{\small 60F10, 60F17, 60G50, 60K05}
\end{abstract}


\bigskip

\section{Introduction}

The concept of a record pertains to instances of achieving the best or worst outcomes within a series of events, holding significance across various domains and relating to all aspects of human life. Notably, records find applications in diverse fields, including finance, environmental studies, physics, and more. For example, in financial contexts, attention is directed towards stock prices reaching apex or nadir points. In the environmental realm, scientists aim to predict the occurrence of peak summer temperatures in various regions, especially as extreme weather events have become increasingly frequent in recent years. Additionally, the study of extreme value theory has enabled scientists to forecast and monitor flood levels. Records also play a pivotal role in physics, as seen in the domain of spin-glasses theory. Consequently, an increasing interest has emerged among researchers in understanding the underlying principles of records, including record generation timing, record numbers, and the extremal values associated with records.

This paper focuses on the number of record occurrences in random walks, including discrete-time and continuous-time. Other papers also refer to the point at which the record is taken as the ladder point. Firstly, we define records at discrete-time random walks. Consider a random walk defined by $S=\left\{S_n=\sum\limits_{i=1}^n X_i; n\ge1, i\ge1\right\}$~($S_0=0$), where $X_i$ represents independent identically distributed random variables. A ``weak record'' is achieved at step $n$ if $S_n=\max\limits_{0\leq k\leq n}S_k$. Similarly, a ``strong record'' occurs if $S_n>\max\limits_{0\leq k\leq n-1}S_k$, indicating a new maximum. For simplicity, related terms mentioned below, like records, ladder points, ladder height, ladder epoch would represent the weak versions. Let $R_n$ represents the number of weak records attained within the first $n$ steps. $R=\{R_n:n\ge0\}$~(where $R_0=1$) stands for the process of record numbers and let $R_{\infty}=\lim_{n\rightarrow\infty}R_n$. The first part of this paper summarizes and augments existing results concerning the weak convergence of $R_n$, the functional limits of $R$, and the distributions of $R_{\infty}$ in scenarios where $S$ drifts towards $-\infty$. Additionally, the paper delves into similar properties about $\sigma$-record, a specific type of record which is more practicable in researches due to the limitation of machine accuracy, the presence of error and noise. Intuitively, it means that we confine records with a threshold $\sigma\ge0$, if the latter record cannot exceed the former $\sigma$-record more than $\sigma$, we would not say that we get a new $\sigma$-record. In addition, we will discuss a sufficient condition for which record numbers, or more appropriately called record time, will have functional limits in continuous-time random walks. The second part of the paper is dedicated to investigating moderate deviations, large deviations, and laws of the iterated logarithm of $R_n$.

Before presenting the results of this paper, it is essential to provide a brief overview of prevailing theories and established conclusions concerning the number of records in discrete-time random walks. For simplicity, we will use ``random walks'' to represent ``discrete-time random walks''. The theoretical foundation of this paper is built upon these insights. Treating the attainment of a record as states within a random walk, it becomes evident that these states are stationary. Consequently, the process $R$ can be regarded as an occupation time process. Drawing on the work of Darling and Kac \cite{Darling-Kac1957} on occupation time in Markov processes, the weak convergence of $R_n$ and the functional limit of $R$ can be explored with suitable scaling. Darling and Kac demonstrated that the normalized occupation time in a stationary Markov process converges weakly to the Mittag-Leffler distribution under a set of uniformity conditions termed D-K conditions. This foundational result has paved the way for subsequent research, as discussed in \cite{bib:pitman1999}. Consequently, the crux of the record problem lies in establishing the D-K conditions. However, despite the significance of Darling and Kac's findings, establishing D-K conditions across a range of scenarios proves challenging.

A pivotal breakthrough in this field originates from Spitzer \cite{bib:fourteen} and \cite{bib:nine}, who derived remarkable identities that could simplify the D-K conditions in random walks. By defining $L_{n,n}$ as the last maximum position within ${S_0,\cdots,S_n}$, Spitzer established a key identity:
\begin{equation}
    \sum_{k=0}^{\infty}\E\left(\e^{-tS_k};L_{k,k}=k\right)y^k = \exp\left\{\sum_{k=1}^{\infty}\frac{y^k}{k}\E\left(\e^{- tS_k};S_k\ge0\right)\right\},
\end{equation}
where $t>0$, $0<y<1$. The left part of the equation serves as a representation of D-K conditions in random walks regarding the number of records. For more concise and understandable proofs of this identity, refer to Karlin and Taylor \cite{bib:two}, page 459. Spitzer's work streamlines D-K conditions by emphasizing aggregate properties of partial sums rather than intricate extreme value problems at each step. After being explored by numerous researchers, Spitzer's condition has emerged as the most central condition in solving the weak convergence problems of $R_n$. The expression of Spitzer's condition is as follows:
\begin{equation}
    \frac{1}{n}\sum_{k=1}^{n}\P\left(S_k\ge0\right) \rightarrow \rho \quad \text{as} \quad n \rightarrow \infty, \quad \rho \in [0,1].
\end{equation}
Following Spitzer's impressive results, Heyde \cite{bib:heyde} and Emery \cite{bib:emery} subsequently demonstrated that if Spitzer's condition holds for some $0<\rho<1$, then $\displaystyle{\sum\limits_{k=0}^{\infty}\P\left(L_{k,k}=k\right)y^k}$ would be a regularly varying function as $y\rightarrow1-$ with parameter $\rho$. Combining Heyde and Emery's findings with D-K conditions, it can be concluded that if Spitzer's condition holds for some $0<\rho<1$, $R_n$ would weakly converge to a Mittag-Leffler distribution with parameter $\rho$, and the regularization should follow the form $n^{-\rho}L(n)$, where $L(n)$ is a slowly varying function. Building on the work of Heyde \cite{bib:heyde}, Emery \cite{bib:emery}, Darling and Kac \cite{Darling-Kac1957}, and Bingham \cite{bib:twlve} \cite{bib:bingham 1972}, Bingham \cite{bib:bingham 1973} established a sufficient and necessary condition for which $R$ would have a non-degenerate ergodic limit, precisely when Spitzer's condition holds for some $0<\rho<1$. Furthermore, Bingham also derived that the limit process is the inverse process of the stable subordinator with index $\rho$, as referred to in \cite{bib:bingham 1973}. In \cite{bib:doney1995}, \cite{doney 1997}, Bertoin and Doney simplified Spitzer's condition further. They proved that $\P\left(S_n\ge0\right)\rightarrow\rho$ is equivalent to $\displaystyle{\frac{1}{n}\sum\limits_{k=1}^{n}\P\left(S_k\ge0\right)\rightarrow\rho}$ in random walks for $\rho\in[0,1]$. At this point, it can be considered that the problems of weak convergence and the functional limit of record numbers are resolved. The core of the remaining questions may lie in investigating the circumstances under which the condition $\P\left(S_n\ge0\right)\rightarrow\rho$ could be satisfied. Undeniably, stable laws have a natural connection with this condition. For example, if the distribution of $\{X_i:~i\ge1\}$, follows a stable distribution or is symmetrical, according to the results, $R_n$ would weakly converge to a truncated normal distribution after regularization, and the functional limit of $R$ could be derived.

In Section 2, we first provide a summary and enhancement of the results concerning the weak convergence of $R_n$ and the functional limits of $R$ under the condition $\P(S_n\ge0)\rightarrow\rho$ for $\rho\in[0,1]$. Subsequently, we extend these findings to $\sigma$-record numbers. Following that, we present a sufficient condition under which record numbers in continuous-time random walks exhibit weak convergence properties and possess non-degenerate functional limits.

Section 3 is dedicated to exploring deviation problems related to $R_n$. Leveraging the renewal property of $R_n$, we employ Cram\'{e}r's method, as outlined in \cite{bib:three}, Theorem 2.2.3, to address large deviations concerning $R_n$, with the result presented in Theorem \ref{ldp}. When Spitzer's condition holds for $0\leq\rho<1$, $R_n$ exhibits weak convergence with a growth rate of $n^{\rho}C_{\rho}(1-1/n)$. Naturally, if we reduce the level of normalization of $R_n$ in large deviation inequality, the moderate deviations about $R_n$ also warrant exploration. In Chen \cite{bib:seven}, moderate deviation problems under D-K conditions have been thoroughly discussed. Gantert and Zeitouni \cite{bib:six} have provided additional results on large deviations and moderate deviations of local time in Markov chains under specific situations. Li and Yao \cite{bib:eleven} have addressed large and moderate deviations in left-continuous or right-continuous random walks with an expectation of $0$. In this paper, we extend the scope to more general random walks. Without imposing any assumptions, for every random walk, regardless of whether it oscillates or drifts, we deduce overall large deviations concerning $R_n$. With the assistance of Chen's results in \cite{bib:seven}, we establish moderate deviations about $R_n$ in Theorem \ref{mdp} under Spitzer's condition. Additionally, by leveraging the peculiarities of random walks, we derive laws of the iterated logarithm of $R_n$ in Theorem \ref{iterated}.

\section{Weak convergence of $R_n$}

We will begin by briefly explaining and introducing the weak convergence results of $R_n$ and the functional limits of $R$ in discrete random walks, where each step length $X_1$ is a non-degenerate random variable, and $S_0$ is considered as the first strong/weak record throughout the paper. As mentioned earlier, for simplicity, the term ``records'' mentioned below will refer to ``weak records''. Let $T_i$ represent the $i$-th ladder epoch in $S$ for $i\ge1$, denoting that $\displaystyle{T_i=\inf\left\{n-\sum\limits_{k=0}^{i-1}T_k;n>\sum_{k=0}^{i-1}T_k, S_n\ge S_j \text{ for all } j=0,1,\dots,n\right\}}$, and set $T_0=0$. Thus, $T_1$ signifies the waiting time to obtain a new record. Referring to \cite{bib:two}, as established, a random walk must fall into one of three states: oscillating, drifting to $+\infty$, or drifting to $-\infty$. This classification depends on the convergence of $\displaystyle{\sum\limits_{k=1}^{\infty}\frac{\P\left(S_k\ge0\right)}{k}}$ and $\displaystyle{\sum\limits_{k=1}^{\infty}\frac{\P\left(S_k<0\right)}{k}}$.
\par
Clearly, when $S$ drifts to $-\infty$, $R_{\infty}$ should be finite with probability $1$. Conversely, when $S$ drifts to $+\infty$ or oscillates, $R_{\infty}$ will almost surely be infinite. Therefore, our attention turns to identifying the weak convergence of $R_n$. Additionally, Bingham \cite{bib:bingham 1973} establishes a sufficient and necessary condition under which $R$ exhibits a non-degenerate functional limit. Remarkably, this condition aligns precisely with Spitzer's condition. The subsequent Theorem \ref{weak convergence theorem} primarily consolidates prior findings. For the part on weak convergence, we provide a concise proof utilizing Darling and Kac's \cite{Darling-Kac1957} Theorem 5:

\begin{theorem}\label{weak convergence theorem}{\rm\bf (Weak convergence and functional limits of record numbers in discrete-time random walks)}
Write $g_{\rho}$ as the Mittag-Leffler distribution with parameter $\rho$, $0\leq \rho \leq1$. There exists some $u(n)>0$ such that $R_n/u(n)$ has a weak convergence limit if and only if $\P\left(S_n\ge0\right)\rightarrow\rho$ with $0\leq\rho\leq1$. In this case, taking $C_{\rho}(1-1/y)=\exp\left[\sum_{k=1}^{\infty}\frac{(1-1/y)^k}{k}\left(\P(S_k\ge0)-\rho\right)\right]$ for $y\ge1$, then $C_{\rho}(1-1/y)$ is a slowly varying function at $+\infty$. Moreover:
\par
(i) When $0<\rho<1$,
$$\frac{R_n}{C_{\rho}(1-1/n)n^{\rho}}\stackrel{d}{\rightarrow} g_{\rho}.$$
Additionally, Bingham notes that the process $R$ would also have a functional limit, which is the inverse process of the stable subordinator with index $\rho$. This condition is also necessary, as detailed in \cite{bib:bingham 1973} and \cite{bib:bingham 1972}.
\par
(ii) When $\rho=0$ with $\displaystyle{\sum\limits_{k=1}^{\infty}\frac{\P\left(S_k\ge0\right)}{k}=+\infty}$,
$$\frac{R_n}{C_0(1-1/n)}\stackrel{d}{\rightarrow} g_0\sim\text{Exponential}\ (1).$$
\par
(iii) When $\rho=0$ with $\displaystyle{\sum\limits_{k=1}^{\infty}\frac{\P\left(S_k\ge0\right)}{k}<+\infty}$, indicating that the random walk $S$ will drift to $-\infty$, then $R_{\infty}$ will follow a geometric distribution with parameter $\displaystyle{\exp\left\{-\sum\limits_{k=1}^{\infty}\frac{\P(S_k\ge0)}{k}\right\}}$, namely
$$\E(R_{\infty})=\exp\left\{\sum\limits_{k=1}^{\infty}\frac{\P(S_k\ge0)}{k}\right\}.$$
\par
(iv) When $\rho=1$,
$$\frac{R_n}{C_1(1-1/n)n}\stackrel{d}{\rightarrow} 1.$$
In this case, $R$ will have a degenerate functional limit. In particular, if $\displaystyle{\sum\limits_{k=1}^{\infty}\frac{\P\left(S_k<0\right)}{k}<\infty}$, indicating that $S$ drifts to $+\infty$, the convergence could be almost surely.
\end{theorem}

\noindent \textbf{Remarks:}
\vspace{0.15cm}

$\bullet$ Bertoin and Doney demonstrated in \cite{bib:doney1995} and \cite{doney 1997} that $\P\left(S_n\ge0\right)\rightarrow\rho$ is equivalent to $\displaystyle{\frac{1}{n}\sum\limits_{k=1}^{n}\P\left(S_k\ge0\right)\rightarrow\rho}$. This equivalence holds for $\rho\in[0,1]$ in random walks as $n\rightarrow\infty$, constituting Spitzer's condition. Therefore, conditions in Theorem \ref{weak convergence theorem} related to $\P\left(S_n\ge0\right)\rightarrow\rho$ could be replaced by Spitzer's condition. Additionally, in (iii) of Theorem \ref{weak convergence theorem}, $\displaystyle{\sum\limits_{k=1}^{\infty}\frac{\P\left(S_k\ge0\right)}{k}<\infty}$ implies $\P\left(S_n\ge0\right)\rightarrow0$, but the reverse may not be true.
\vspace{0.2cm}

$\bullet$ Referring to Ros\'{e}n \cite{bib:rosen}, Theorem 1, in random walks, if $X_1$ follows a non-degenerate distribution, then $\P(S_n=0)\leq O(1/\sqrt{n})$. Thus, the distinction between $\P(S_n>0)$ and $\P(S_n\ge0)$ can be disregarded in the conditions of Theorem \ref{weak convergence theorem}. Additionally, if we replace $C_{\rho}$ with $\hat{C}{\rho}$, where $\displaystyle{\hat{C}{\rho}(y)=\exp\left\{\sum\limits_{k=1}^{\infty}\frac{y^k}{k}\left(\P(S_k>0)-\rho\right)\right\}}$ for $0<y<1$, Theorem \ref{weak convergence theorem} would hold for strong record numbers. It can be observed that the difference between weak record numbers and strong record numbers is just a constant multiple.
\vspace{0.2cm}

By leveraging the Central Limit Theorem and Berry-Ess\'{e}en's Theorem, we can derive specific examples:
\begin{example}\label{berry-esseen}
If the distribution of $X_1$ has finite third moment with $\mu=\E(X_1)=0$~(such as when $X_1$ is bounded), there exists a positive constant $C_{1/2}^1$ such that
$$\frac{R_n}{C^1_{1/2}n^{1/2}}\stackrel{d}{\rightarrow}g_{1/2}.$$
In the case where the distribution of $X_1$ is continuous and symmetrical around $0$, then
$$\frac{R_n}{n^{1/2}}\rightarrow g_{1/2}.$$
\end{example}

\begin{example}[An example in Paulauskas \cite{bib:stable law}]
If $X_1$ follows the density function
$$f(x)=\frac{1-\cos x}{\pi x^2},$$
and with reference to Paulauskas \cite{bib:stable law} and Christoph \cite{bib:thirteen}, it can be established that $X_1$ belongs to the domain of attraction of a strict stable distribution $St(1,1/2)$. Moreover, if the rate of convergence is sufficiently fast to satisfy $\displaystyle{\sum\limits_{k=1}^{\infty}\frac{1}{k}|\P(S_k\ge0)-\rho|<+\infty}$, then there exists a positive constant $C^2_{1/2}$ such that
$$\frac{R_n}{C^2_{1/2}n^{1/2}}\stackrel{d}{\rightarrow}g_{1/2}.$$
\end{example}

\begin{example}[Li and Yao \cite{bib:eleven}]\label{Yuqiang}
Assuming $X_1$ is left-continuous with $\P(X_1=k)=p_k$ for $k=-1,0,1,\dots$ and $\E(X_1)=0$, we define
$$\varphi(s)=\sum_{k=-1}^{\infty}p_{k}s^{k+1}$$
for $0\leq s\leq1$. If
$$\varphi(s)=s+\frac{\gamma}{1+\beta}(1-s)^{1+\beta}$$
for some $\beta\in(0,1)$ and $\gamma\in(0,1)$, then $X_1$ falls within the domain of attraction of a strict stable distribution $St(1+\beta,\frac{1}{1+\beta})$.  Consequently, Theorem \ref{weak convergence theorem} holds with $\rho=1/(1+\beta)$ and $C_{\rho}=\left(\frac{1+\beta}{\gamma}\right)^{\beta/(1+\beta)}$.
\end{example}

\noindent \textbf{Remarks:}
\vspace{0.15cm}

$\bullet$ In Li and Yao \cite{bib:eleven}, they introduced assumption (H) as a condition ensuring the weak convergence of $R_n$ when $X_1$ is a left-continuous or right-continuous integral-valued random variable. Specifically, under the assumption (H) in \cite{bib:eleven}, for the case of $X_1$ being left-continuous,  it implies that
\begin{equation}\label{yuqiang assumption h}
\exp\left[\sum_{k=1}^{\infty}\frac{y^k}{k}\P(S_k<0)\right]\sim \frac{p_{-1}}{c(1-y)^{\alpha}}\end{equation}
as $y\rightarrow1-$, where $c>0$ is a constant, $\alpha\in(0,1)$, and $p_{-1}=\P(X_1=-1)$. Refer to \cite{bib:eleven} for the definitions of these notations. Drawing on earlier discussions and Lamperti \cite{bib:lamperti}, equation (\ref{yuqiang assumption h}) should imply $\P(S_n\ge0)\rightarrow1-\alpha$.
\vspace{0.2cm}

Furthermore, if $N_n$ and $N^+_n$ represent the counts of $S_i \ge 0$ and $S_i > 0$ within $\{S_1,\dots,S_n\}$ respectively, it is known that Spitzer established a notable arc-sine law concerning $N_n$ in \cite{bib:nine} under Spitzer's condition. Drawing inspiration from \cite{bib:two}, page 454, and leveraging equivalence principles within random walks, along with the non-increasing nature of $\P(N_n=n)$ and $\P(N_n^+=n)$, combined with Tauberian's theorem, more precise estimations about $\P(L_{n,n}=n)$ and $\P(N_n=n)$ can be derived when $n$ is sufficiently large. Here, $L_{n,0}$ denotes the position of the first maximum among $\{S_0,\cdots,S_n\}$, and attentive readers will quickly discern its connection to strong records.

\begin{corollary}\label{limit of Ln,n}
If $\P\left(S_n\ge0\right)\rightarrow\rho$ as $n\rightarrow\infty$ for $0<\rho\leq1$, then
$$\P(L_{n,n}=n)=\P(N_n=n)\sim\frac{1}{\Gamma(\rho)}n^{\rho-1}C_{\rho}(1-1/n),$$ and
$$\P(L_{n,0}=n)=\P(N_n^+=n)\sim\frac{1}{\Gamma(\rho)}n^{\rho-1}\hat{C}_{\rho}(1-1/n).$$
In particular, when $\displaystyle{\sum\limits_{k=1}^{\infty}\frac{1}{k}\P\left(S_k<0\right)<+\infty}$,
$$\lim_{n\rightarrow\infty}\P(L_{n,n}=n)=\lim_{n\rightarrow\infty}\P(N_n=n)=\exp\left[-\sum_{k=1}^{\infty}\frac{1}{k}\P\left(S_k<0\right)\right],$$ and
$$\lim_{n\rightarrow\infty}\P(L_{n,0}=n)=\lim_{n\rightarrow\infty}\P(N_n^+=n)=\exp\left[-\sum_{k=1}^{\infty}\frac{1}{k}\P\left(S_k\leq0\right)\right].$$
\end{corollary}

In cases where the support of $X_1$ is integer and $X_1$ is right-continuous, indicating that $\P(X_1=1)>0$ while $\P(X_1=k)=0$ for $k=2,3,\dots$, there exists an intuitive connection between the number of strong records and the maximum value. Let $M_n$ be defined as the maximum value among $\{S_0,S_1,\dots,S_n\}$, and consider  $\displaystyle{M_{\infty}=\lim\limits_{n\rightarrow\infty}M_n}$.

\begin{corollary}
Assuming $X_1$ is right-continuous, if $\P\left(S_n\ge0\right)\rightarrow\rho$ holds with $0<\rho\leq1$ or $\P\left(S_n\ge0\right)\rightarrow0$ with $\sum_{k=1}^{\infty}\frac{\P\left(S_k\ge0\right)}{k}=\infty$, then $M_n/[\hat{C}_{\rho}(1-1/n)n^{\rho}]\stackrel{d}{\rightarrow}g_{\rho}$.
\par

If $\sum_{k=1}^{\infty}\frac{\P(S_k\ge0)}{k}<\infty$, then $M_{\infty}<+\infty$ with probability $1$, and it follows a geometric distribution with an expectation of $\displaystyle{\exp\left\{\sum\limits_{k=1}^{\infty}\frac{\P(S_k>0)}{k}\right\}}$.
\end{corollary}

At times, due to considerations of machine testing accuracy or experimental design requirements, subtle variations in records may be challenging to detect or may not significantly impact the experiment, qualifying as noise. Consequently, researchers find it necessary to establish a threshold for what qualifies as a new record. This concept is known as $\sigma$-records. In mathematical terms, if $S_k$ is a $\sigma$-record, for $j>k$, $S_j$ becomes the first weak $\sigma$-record following $S_k$ provided that $S_j - S_k \ge \sigma$. A similar distinction exists between weak $\sigma$-records and strong $\sigma$-records. Without special descriptions, below, ``$\sigma$-record'' would refer to ``weak $\sigma$-record''. Similarly, regarding $S_0$ as the first weak $\sigma$-record point and defining $R^{\sigma}_n$ as the number of weak $\sigma$-records during the first $n$ steps (with $R^{\sigma}_0=1$), taking $R^{\sigma}=\{R^{\sigma}_n:n\ge0,R^{\sigma}_0=1\}$, we will delve into the matter of the weak convergence of $R^{\sigma}_n$ and the functional limit problem concerning $R^{\sigma}$, employing a similar approach. Let $Z_i=S_{T_i}-S_{T_{i-1}}$, $i\ge1$, $H_n = \sum_{k=0}^{n}Z_k$ and $H_0=0$. Then, defining
$$V(x) = \sum_{n=0}^{\infty}\P(H_n\leq x),$$
we can consider $V(x)$ as a renewal function associated with $Z_1$.

\begin{theorem}\label{sigma-record}{\rm\bf(Weak convergences and functional limits of $\sigma$-record numbers in discrete-time random walks)}
Assume $\sigma$ is a continuity point of $V(x)$:
\par
(i) If $\P\left(S_n\ge0\right)\rightarrow\rho$ with $0<\rho<1$,
$$\frac{V(\sigma)R^{\sigma}_n}{C_{\rho}(1-1/n)n^{\rho}}\stackrel{d}{\rightarrow} g_{\rho},$$
and similarly, the process $R^{\sigma}$ will have a functional limit that is the inverse process of the stable subordinator with index $\rho$.
\par
(ii) If $\P\left(S_n\ge0\right)\rightarrow0$ with $\displaystyle{\sum\limits_{k=1}^{\infty}\frac{\P\left(S_k\ge0\right)}{k}=\infty}$,
$$\frac{V(\sigma)R^{\sigma}_n}{C_{\rho}(1-1/n)}\stackrel{d}{\rightarrow} g_0.$$
If $\P\left(S_n\ge0\right)\rightarrow1$, $R^{\sigma}_n$ would have a degenerate weak convergence limit with normalization $\displaystyle{\frac{V(\sigma)}{nC_1(1-1/n)}}$.
\par
(iii) If $S$ drifts to $-\infty$, $\displaystyle{R_{\infty}^{\sigma}=\lim\limits_{n\rightarrow\infty}R_n^{\sigma}}$ will follow a geometric distribution with parameter $\displaystyle{V(\sigma)\exp\left\{-\sum\limits_{k=1}^{\infty}\frac{\P(S_k\ge0)}{k}\right\}}$.
\end{theorem}
\noindent \textbf{Remarks:}
\vspace{0.15cm}

$\bullet$ $V(x)$ has various equivalent forms; refer to \cite{doney 1994} for detailed explanations. Additionally, according to \cite{bib:bingham 1973}, the Laplace-Stieltjes transform of $V$ is given by $\displaystyle{\exp\left\{\sum\limits_{k=1}^{\infty}\frac{1}{k}\E\left(\e^{-\lambda S_k};S_k\ge0\right)\right\}}$, where the independent variable is denoted by $\lambda>0$.
\vspace{0.2cm}

Consider a continuous-time random walk $\widetilde{S}=\left\{\widetilde{S}_n;n\ge0\right\}$ with i.i.d. step length  $X_i,$ and i.i.d. waiting time $Y_i$ for each step, where $\displaystyle{\widetilde{S}_t=\sum\limits_{i=1}^{\max\left\{k:\sum_{j=1}^{k}Y_j\leq t\right\}}X_i}$ for $t\in \mathsf{R}^+$, and $\widetilde{S}_0=0$. Similarly, the record numbers, or record times obtained before time $t$ in continuous-time random walks, are denoted as $\widetilde{R}_t$. Let $\widetilde{R}=\{\widetilde{R}_t;t\ge0\}$ with $\widetilde{R}_0=1$. Additionally, we assume that $X_i$ and $Y_i$ are independent for $i\ge 1$ in the following theorem, which provides a sufficient condition for the record time $\widetilde{R}_t$ to have weak convergence limits and $\widetilde{R}$ to have functional limits. We still denote $\displaystyle{S_n=\sum\limits_{k=1}^nX_k}$.
\begin{theorem}\label{ctrw records weak convergence}{\rm\bf(Weak convergence and functional limits of record numbers in continuous-time random walks)}
If $\P(S_n\ge0)\rightarrow\rho$ for $0<\rho\leq1$, $Y_1$ is in the domain of attraction of a stable distribution with parameter $\alpha\in(0,1)$, meaning that the cumulative distribution function of $Y_1$, denoted as $G(x)$, satisfies
$$1-G(x)\sim \frac{1}{\Gamma(1-\alpha)}x^{-\alpha}L_1(x),$$
where $L_1(x)$ is a slowly varying function at $+\infty$, then
$\widetilde{R}$ has a functional limit that is the inverse process of the stable subordinator with index $\alpha\rho$.

\end{theorem}

\section{Large, Moderate Deviations and Laws of the Iterated Logarithm}

We have previously established that $R$ can be conceptualized as a type of renewal process, making theories of large deviations applicable to some extent. As referenced in \cite{bib:two} on page 465, we are already acquainted with the fact that for $|y| < 1$,
$$\E(y^{T_1}) = 1 - \exp\left\{-\sum_{k=1}^{\infty}\frac{y^k}{k}\P\left(S_k \ge 0\right)\right\}.$$
Through elementary deductions, we can further deduce that $\E(T_1) = \infty$ if and only if $\displaystyle{\sum\limits_{k=1}^{\infty}\frac{\P(S_k < 0)}{k} = +\infty}$, and $T_1$ can be regarded as the waiting time in $R$. On one hand, if $\E(T_1) < +\infty$, the large deviations about $R_n$ are just formal problems of large deviations in renewal processes. On the other hand, in scenarios where $\displaystyle{\sum\limits_{k=1}^{\infty}\frac{\P(S_k \ge 0)}{k} < +\infty}$, as per Theorem \ref{weak convergence theorem}, we establish that $R_{\infty}$ will be finite with probability $1$. Therefore, discussing large deviations in this case becomes meaningless.
\par
Define $\Lambda(\lambda)=\log\E(\e^{\lambda T_1})$ as $T_1$'s cumulant generation function, and $\Lambda^*({y})=\sup_{\lambda\leq0}\left(\lambda y-\Lambda(\lambda)\right)$. We present the large deviations about $R_n$ as follows.
\begin{theorem}[Large deviations]\label{ldp}
If the random walk $S$ oscillates, the distributions of $R_n/n$ satisfy a large deviation principle: for any $0 < y \leq 1$,
$$\lim_{n\rightarrow\infty}\frac{1}{n}\log\P\left(R_n\ge ny\right)=-y\Lambda^*(\frac{1}{y}),$$
where
$$\Lambda^*(y)=\begin{cases}
(\Lambda')^{-1}(y)\cdot y-\Lambda\left[(\Lambda')^{-1}(y)\right],&y>1,\\
-\log\P(X_1\ge0),&y=1,\\
+\infty,&y<1,
\end{cases}$$
where $\Lambda'(y)$ is the derivative of $\Lambda(y)$, and $(\Lambda')^{-1}(y)$ is the inverse of $\Lambda'(y)$ under $y<0$.
\par
When $S$ drifts to $+\infty$, we already know that $\displaystyle{\frac{R_n}{n}\stackrel{\text{a.s.}}{\rightarrow}\frac{1}{\E(T_1)}}$. Therefore, for $1/\E(T_1)<y\leq1$,
$$\lim_{n\rightarrow\infty}\frac{1}{n}\log\P\left(R_n\ge ny\right)=-y\Lambda^*(\frac{1}{y}),$$
$$\Lambda^*(y)=\begin{cases}
(\Lambda')^{-1}(y)\cdot y-\Lambda\left[(\Lambda')^{-1}(y)\right],&\E(T_1)>y>1,\\
-\log\P(X_1\ge0),&y=1,\\
+\infty,&y<1,
\end{cases}$$
where $\displaystyle{\E(T_1)=\exp\left\{\sum\limits_{k=1}^{\infty}\frac{\P(S_k<0)}{k}\right\}<+\infty}$.
\end{theorem}

Given that $R_n$ has a weak convergence limit, and the regularization parameter is $n^{-\rho}C^{-1}_{\rho}(1-1/n)$ under Spitzer's condition, we can also explore the phenomenon of moderate deviations of $R_n$. To some extent, the moderate deviation problems related to the number of records can be reframed in the context of the occupation time in a Markov chain, linked with Green's function. Numerous results already exist in this domain, as seen in \cite{bib:six} and \cite{bib:seven}.

Before presenting our findings, we need to introduce a new symbol $a(n)$, which serves to regulate the rate of change of small probability events within the context of moderate deviation problems. Additionally, specific constraints must be imposed on $a(n)$ depending on the speed of the regularization parameter of $R_n$ in the context of weak convergence.

\begin{theorem}[Moderate deviations]\label{mdp}
If $\P\left(S_n\ge0\right)\rightarrow\rho$ for $0<\rho<1$ or $\P\left(S_n\ge0\right)\rightarrow0\ (\rho=0)$ with $\sum_{k=1}^{\infty}\frac{\P(S_k\ge0)}{k}=\infty$. Assuming that $a(n)\rightarrow\infty$ and
$\frac{a(n)}{n^{1-\rho}}\rightarrow0$
as $n\rightarrow\infty$, the distributions of $\frac{R_n}{C_{\rho}\left(1-\frac{a(n)^{1/(1-\rho)}}{n}\right)n^{\rho}a(n)}$ satisfy a moderate deviation principle: for $y>0$,
$$\lim_{n\rightarrow\infty}\frac{1}{a(n)^{1/(1-\rho)}}\log\left[\P\left(\frac{R_n}{C_{\rho}\left(1-\frac{a(n)^{1/(1-\rho)}}{n}\right)n^{\rho}a(n)}\ge y\right)\right]=-(1-\rho)\left(\rho^{\rho}y\right)^{1/(1-\rho)}.$$
\end{theorem}

We also present our version of laws of the iterated logarithm regarding record numbers.

\begin{theorem}[Laws of the iterated logarithm]\label{iterated}
If $\P\left(S_n\ge0\right)\rightarrow\rho$ for $0<\rho<1$ or $\P\left(S_n\ge0\right)\rightarrow0$ with $\sum_{k=1}^{\infty}\frac{\P(S_k\ge0)}{k}=\infty$, for any $(x_1,x_2)\in \mathsf{R^2}$, define $R^{(x_1,x_2)}_n$ as the number of weak records satisfying that the records are no smaller than $x_2$ during the first $n$ steps, and the random walk $S$ starts at $x_1$. Then,
$${\lim\sup}_{n\rightarrow\infty}\frac{R^{(x_1,x_2)}_n}{f\left(n/\log\log f(n)\right)\log\log f(n)}=\frac{\Gamma(\rho+1)}{\rho^{\rho}(1-\rho)^{1-\rho}}\ \ a.s.,$$
where $f(n)=\frac{n^{\rho}C_{\rho}(1-1/n)}{\Gamma(\rho+1)}$.
\end{theorem}

\noindent \textbf{Remarks:}
\vspace{0.15cm}

$\bullet$ The large, moderate deviations and laws of the iterated logarithm also hold for strong record numbers $\hat{R}_n$ with some minor adjustments.

\section{Proof of Weak Convergence and Functional Limits}

\noindent\textbf{Proof of Theorem \ref{weak convergence theorem}.}~The main part of Theorem \ref{weak convergence theorem} originates from Theorem 1 in Bingham \cite{bib:bingham 1973}. Some augmentations in Theorem \ref{weak convergence theorem} can be straightforwardly verified using Darling-Kac's occupation time theory in \cite{Darling-Kac1957}, thus we provide a concise explanation. Here, we offer an alternative perspective to prove (iii) and provide additional details regarding the almost sure part in (iv).
\par
A brief explanation regarding the weak convergence part in (i) and (ii):
\par
Define $\Theta=\{(x_1,x_2):x_1=x_2 \in \mathsf{Z},x_1\ge0\}$ and a new bivariate Markov chain $(S,M)=\{(S_n,M_n);n\ge0\}$. For $0<y<1$, $0\leq \rho\leq1$, using the identities in Karlin's \cite{bib:two}, page 459,
$$\begin{aligned}
\sum_{n=0}^{\infty}\P\left((S_{n+1},M_{n+1})\in\Theta|(S_0,M_0)=(0,0)\right)y^n&=\sum_{n=1}^{\infty}\P\left(L_{n,n}=n\right)y^{n-1}\\
&=\frac{\exp\left(\sum_{k=1}^{\infty}\frac{y^k}{k}\P\left(S_k\ge0\right)\right)-1}{y}\\
&= \frac{C_{\rho}(y)(1-y)^{-\rho}-1}{y},
\end{aligned}$$
It is known that $C_{\rho}(1-1/n)$ is a slowly varying function as $n$ approaches $+\infty$ if and only if $\P(S_n\ge0)\rightarrow\rho$, as established by Emery \cite{bib:emery} or Bingham \cite{bib:bingham 1973}. Additionally, without a doubt, the occupation time of $\Theta$ in $(S,M)$ during the first $n$ steps equals $R_n$. For any $(x,x)\in\Theta$, $p((x,x),\Theta,y)=p((0,0),\Theta,y)$ holds due to the renewal property of $R_n$, signifying uniformity in D-K conditions. Therefore, by applying Darling-Kac's occupation time theory, we obtain the weak convergence results.

\par

A proof of (iii):
\par
When $\sum_{k=1}^{\infty}\frac{\P(S_k\ge0)}{k}<\infty$, from Spitzer's \cite{bib:nine}, Theorem 4.1, $M_{\infty}$ would be finite with probability 1 and
$$\E(y^{M_{\infty}})=\prod_{k=1}^{\infty}\exp\left[\frac{\psi_k(y)-1}{k}\right],$$
for $0<y<1$ where $\psi_k(y)=\E(y^{S_k^+})$, $S_k^+=\max\{S_k,0\}$. Besides,
$$\E\left(\left.y^{Z_1}\right|Z_1<\infty\right)=\left\{1-\exp\left[-\sum_{k=1}^{\infty}\frac{1}{k}\E(y^{S_k};S_k\ge0)\right]\right\}/\P(Z_1<\infty),$$
and here $\{Z_1<\infty\}$ represents that $S$ would have next record, which is also equivalent with $\{T_1<\infty\}$. In addition,
$$\P(Z_1=\infty)=\P(T_1=\infty)=\exp\left[-\sum_{k=1}^{\infty}\frac{\P(S_k\ge0)}{k}\right].$$
Since $Z_i$, $i\ge1$ are i.i.d. from the Markov property of $S$, $\sum_{i=1}^{R_{\infty}}Z_i=M_{\infty}$, using the independence of $Z_i$ and $R_{\infty}$, referring to \cite{bib:two} or \cite{bib:regluar book},
$$\begin{aligned}\E(y^{M_{\infty}})&=\sum_{n=1}^{\infty}\P(R_{\infty}=n)\E\left(\left.y^{Z_1}\right|Z_1<\infty\right)^{n-1}\\
&=\sum_{n=1}^{\infty}\P(R_{\infty}=n)(1-\P(Z_1=\infty))^{1-n}\left\{1-\exp\left[-\sum_{k=1}^{\infty}\frac{\psi_k(y)-\P(S_k<0)}{k}\right]\right\}^{n-1},\end{aligned}$$
then we could deduce that $R_{\infty}$ has the geometric distribution with $\E(R_{\infty})=\exp\left[\sum_{k=1}^{\infty}\frac{\P(S_k\ge0)}{k}\right]$.
\par
Using renewal theory to explicate a.s. convergence type if $\sum_{k=1}^{\infty}\frac{\P(S_k<0)}{k}<\infty$:
\par
Referring to Karlin's \cite{bib:two} page 465, Theorem 4.1 and page 481,
$$\E(T_1)=\exp\left[\sum_{k=1}^{\infty}\frac{\P(S_k<0)}{k}\right].$$
It can be verified that the process of producing a new record is a renewal process with waiting time $T_1$. With $E(T_1)<\infty$, using S.L.L.N. in renewal process we get a.s.
$$\lim_{n\rightarrow\infty}\frac{R_n}{n}=\frac{1}{\E(T_1)}.$$
Since that $g_1$ is the distribution of degenerate random variable on 1, we get the result.

\bigskip

\noindent\textbf{Proof of Theorem \ref{sigma-record}.}~(i) The proof of Theorem \ref{sigma-record} is exactly similar with Proof of Theorem \ref{weak convergence theorem} by considering Emery's \cite{bib:emery}, Bingham's Theorem 2 in \cite{bib:bingham 1973} or Theorem 8.9.12 in \cite{bib:regluar book}. Because $\P(M_n<\sigma)=\P(T^{\sigma}_1> n)$ if we take $T^{\sigma}_1$ as the waiting time of $S$ getting a new $\sigma$-record, according to Bingham's Theorem 2 in \cite{bib:bingham 1973},
$$\P(T^{\sigma}_1> n)\sim \frac{V(\sigma)}{n^{\rho}C_{\rho}(1-1/n)\Gamma(1-\rho)}.$$
Then use Bingham's Theorem 3a in \cite{bib:bingham 1972} we have proved the result.
\par
(ii) It's just extreme situations of (i) except that we won't have non-degenerate functional limits.
\par
(iii) Under the condition, we just need to verify that $\P(T_1^{\sigma}=+\infty)=\frac{V(\sigma)}{C_0}$. Actually, when $\rho=0$ and $C_0<+\infty$,
$$\P(T_1^{\sigma}=+\infty)=\lim_{n\rightarrow\infty}\P(T^{\sigma}_1> n)=\frac{V(\sigma)}{C_0}.$$

\noindent\textbf{Proof of Theorem \ref{ctrw records weak convergence}.}~For a regenerative phenomenon $Z(u,w)$, which is an indicator process, meaning that $\{Z(u,w)=1\}$ represents at time $u$ the event happens (here is our records process). Take $H(t,w)=\int_0^tZ(u,w)\mathrm{d}u$ ($H(t)$ is the occupation time before time $t$), and take $T(v,w)=\sup\{t:H(t,w)\le v\}$, which is the inverse process of $H$. Denote $r(s)$ as the Laplace transform of $p(t)$ where $p(t)=\P(Z(t,w)=1)$,
$$r(s)=\int_0^{+\infty}\e^{-st}p(t)\mathrm{d}t,$$
Under the assumption that $\lim_{t\rightarrow0}p(t)=1$, which is equivalent with our convention that $\widetilde{S}_0$ is the first record, in \cite{bib:bingham 1972}, Bingham have the result that $H$ has ergodic limit if and only if
$$r(s)\sim s^{-\alpha}L(1/s),$$
$s\rightarrow0$, where $L$ is a slowly varing function and $0\leq \alpha <1$.
\par
Also, Kingman has proved in \cite{bib:kingman} that
$$\E(\exp(-sT(v,w)))=\exp(-v/r(s)).$$
Therefore, we could find
$$r(s)\sim s^{-\alpha}L(1/s)$$
is equal with
$$1-\E(\exp(-sT(1,w)))\sim s^{\alpha}L^{-1}(1/s)$$
as $s\rightarrow0$.
\par
We could realize that if we take $\widetilde{T}_1$ as the waiting time of generating a new record at a continuous-time random walk, then
$$\widetilde{T}_1=T(1,w).$$
Take Laplace transform of $Y_1$ as $\phi(t)$, namely
$$\phi(t)=\E\left(\e^{-tY_1}\right)$$
for $t>0$. By Tauberian theorem we get
$$1-\phi(t)\sim t^{\alpha}L_1(1/t)$$
when $t\rightarrow0$. In this case, the Laplace transform of $\widetilde{T}_1$, taking as $f(t)=\E\left(\e^{-t\widetilde{T}_1}\right)$, should be
$$f(t)=1-\exp\left\{-\sum_{k=1}^{+\infty}\frac{\phi^k(t)}{k}\P(S_k\ge0)\right\},$$
this is because when the random walk is discrete, $f(t)=1-\exp\left\{-\sum_{k=1}^{+\infty}\frac{t^k}{k}\P(S_k\ge0)\right\}$ refering to Karlin's \cite{bib:two}.
Under our conditions, as $t\rightarrow0$,
$$1-f(t)\sim t^{\alpha\rho}\hat{L}(t),$$
where $\hat{L}(t)$ is a slowly varying function at $+\infty$. Therefore, by Bingham's result, $\widetilde{R}_t/[\hat{L}(1-1/t)t^{\alpha\rho}]$ weak converge to $g_{\alpha\rho}$ and $\widetilde{R}$ has functional limit, the normalization is just $\hat{L}(1-1/t)t^{\alpha\rho}$.

\section{Proof of Large, Moderate Deviation Principles and Laws of the Iterated Logarithm}

\noindent\textbf{Proof of Theorem \ref{ldp}.}~For any $0< y\leq1$, it can be verified that $W_i,i\ge0$ are stopping times where $W_i=\sum_{k=0}^{i}T_k$ and from the strong Markov property of $S$,
$$\P\left(\frac{R_n}{n}\ge y\right)=\P\left(W_{\lceil yn \rceil-1}\leq n\right).$$
Besides, $T_i, i\ge 1$ are i.i.d. random variables having the same distribution with $W_1=T_1$. By Cram\'er's theorem in \cite{bib:three} we have
$$\lim_{n\rightarrow\infty}\frac{1}{\lceil yn \rceil-1}\log\P\left(\sum_{i=1}^{\lceil yn \rceil-1}T_i\leq n\right)=-\Lambda^*(\frac{1}{y}),$$
where the right part of the equality comes from the non-increasing of $\Lambda^*$, then the result could be obtained directly.
\par
We'll illustrate the form of $\Lambda^*(\lambda)$ below briefly:
\par
We have already know that for $\lambda\leq0$,
$$\Lambda(\lambda)=\log\left[1-\exp\left\{-\sum_{k=1}^{\infty}\frac{\e^{\lambda k}}{k}\P(S_k\ge0)\right\}\right],$$
so
$$\Lambda'(\lambda)=\frac{\left[\sum_{k=1}^{\infty}\e^{\lambda k}\P(S_k\ge0)\right]\exp\left\{-\sum_{k=1}^{\infty}\frac{\e^{\lambda k}}{k}\P(S_k\ge0)\right\}}{1-\exp\left\{-\sum_{k=1}^{\infty}\frac{\e^{\lambda k}}{k}\P(S_k\ge0)\right\}}.$$
As $\lambda\rightarrow-\infty,$ using L'Hospital's rule and $\P(S_1\ge0)>0$ due to that $S$ oscillates or drifts to $+\infty$,
$$\begin{aligned}
\lim_{\lambda\rightarrow-\infty}\Lambda'(\lambda)&=\lim_{\lambda\rightarrow-\infty}\frac{\left[\sum_{k=1}^{\infty}k\e^{\lambda k}\P(S_k\ge0)-\left(\sum_{k=1}^{\infty}\e^{\lambda k}\P(S_k\ge0)\right)^2\right]\exp\left\{-\sum_{k=1}^{\infty}\frac{\e^{\lambda k}}{k}\P(S_k\ge0)\right\}}{\left[\sum_{k=1}^{\infty}\e^{\lambda k}\P(S_k\ge0)\right]\exp\left\{-\sum_{k=1}^{\infty}\frac{\e^{\lambda k}}{k}\P(S_k\ge0)\right\}}\\
&=\lim_{\lambda\rightarrow-\infty}\frac{\left[\sum_{k=1}^{\infty}k\e^{\lambda k}\P(S_k\ge0)-\left(\sum_{k=1}^{\infty}\e^{\lambda k}\P(S_k\ge0)\right)^2\right]}{\left[\sum_{k=1}^{\infty}\e^{\lambda k}\P(S_k\ge0)\right]}\\
&=1.
\end{aligned}$$
If $S$ oscillates, then $\sum_{k=1}^{\infty}\frac{1}{k}\P(S_k<0)=+\infty$, which means $\E(T_1)=+\infty$, and
$$\lim_{\lambda\rightarrow0-}\Lambda'(\lambda)=\lim_{\lambda\rightarrow0-}\frac{\E\left(T_1\e^{\lambda T_1}\right)}{\E\left(\e^{\lambda T_1}\right)}=\E(T_1)=+\infty.$$
Besides,
$$\Lambda''(\lambda)=\frac{\E\left(T_1^2\e^{\lambda T_1}\right)\E\left(\e^{\lambda T_1}\right)-\E^2\left(T_1\e^{\lambda T_1}\right)}{\E^2\left(\e^{\lambda T_1}\right)}\ge0$$
by H\"older's inequality, so $\Lambda'(\lambda)$ should be non-decreasing from 1 to $+\infty$ when $\lambda$ goes to $0$ from $-\infty$. Hence, if $y>1$, $\Lambda^*(y)=(\Lambda')^{-1}(y)\cdot y-\Lambda\left[(\Lambda')^{-1}(y)\right]$; if $y\leq1$,
$$\begin{aligned}
\Lambda^*(y)&=\sup_{\lambda\leq0}\left(\lambda y-\Lambda(\lambda)\right)\\
&=\lim_{\lambda\rightarrow-\infty}\left(\lambda y-\Lambda(\lambda)\right)\\
&=\lim_{\lambda\rightarrow-\infty}\log\left[\frac{\e^{\lambda y}}{1-\exp\left\{-\sum_{k=1}^{\infty}\frac{\e^{\lambda k}}{k}\P(S_k\ge0)\right\}}\right]\\
&=\lim_{\lambda\rightarrow-\infty}\log\left[\frac{y\e^{\lambda y}}{\left[\sum_{k=1}^{\infty}\e^{\lambda k}\P(S_k\ge0)\right]\exp\left\{-\sum_{k=1}^{\infty}\frac{\e^{\lambda k}}{k}\P(S_k\ge0)\right\}}\right]\\
&=\begin{cases}
-\log\P(S_1\ge0),&y=1\\
+\infty.&y<1
\end{cases}
\end{aligned}$$
If $S$ drifts to $+\infty$, then $1\leq\E(T_1)<\infty$, the proof is similar.

On moderate deviation problems, Chen Xia \cite{bib:seven} has given a useful result under D-K conditions:
\begin{lemma}[Chen Xia \cite{bib:seven}]\label{chen}
\ \par
If there exists a function $f(t)$ which is regular varing at infinity with exponent $0\leq p<1$, namely $\lim_{t\rightarrow\infty}\frac{f(\lambda t)}{f(t)}=\lambda^p$ for any $\lambda>0$, and $b_n\rightarrow\infty$, $b_n/n\rightarrow0$ when $n\rightarrow\infty$, then if
$$\lim_{n\rightarrow\infty}\frac{1}{f(n)}\sum_{k=1}^{n}\P\left(S_k\in E\right)=1$$
uniformly in $E$, we have
$$\lim_{n\rightarrow\infty}\frac{1}{b_n}\log\P\left[\sum_{k=1}^{n}I_E(S_k)\ge yf(n/b_n)b_n\right]=-(1-p)\left(\frac{p^p\lambda}{\Gamma(p+1)}\right)^{(1-p)^{-1}}.$$
$I_E$ represents the indicator function of set $E$.
\end{lemma}

\noindent\textbf{Proof of Theorem \ref{mdp}.}~By Proof of Theorem \ref{weak convergence theorem} and Tauberian's theorem we have
\begin{equation}\label{tauberian C1}\sum_{k=1}^{n}\P((S_{k},M_{k})\in\Theta|(S_0,M_0)=(0,0))\sim\frac{n^{\rho}C_{\rho}(1-1/n)}{\Gamma(\rho+1)}
\end{equation}
as $n\rightarrow\infty$. Combining with Chen Xia's result in Lemma \ref{chen}, taking $b_n=a(n)^{\frac{1}{1-\rho}}$, and from (\ref{tauberian C1}), we can take $f(n)=\frac{n^{\rho}C_{\rho}(1-1/n)}{\Gamma(\rho+1)}$ which is a regularly varied function with exponent $\rho$, then we get the result.

\noindent\textbf{Proof of Theorem \ref{iterated}.}~Referring to Theorem 3 in Chen Xia's \cite{bib:seven}, and define
$$R(\Theta)=\left\{(x_1,x_2)\in R^2;\lim_{n\rightarrow\infty}\frac{\Gamma(\rho+1)}{n^{\rho}C_{\rho}(1-1/n)}\sum_{k=1}^{n}\P_{x_1}(L_{k,k}=k,S_k\ge x_2)=1,x_2\ge x_1\right\},$$
where the suffix $x_1$ represents that $S_0=x_1$. When $x_2\leq x_1$, the result is trivial just by using Theorem 3 in Chen Xia's \cite{bib:seven}.
For the remained part, we only need to explain that for every $x_2>0$,
$$\lim_{n\rightarrow\infty}\frac{\Gamma(\rho+1)}{n^{\rho}C_{\rho}(1-1/n)}\sum_{k=1}^{n}\P(L_{n,n}=n,S_n\ge x_2)=1.$$
Denote $\tau_{x_2}$ as the first passage time to $[x_2,+\infty)$ from initial point 0, then for any $0<y<1$,
$$\begin{aligned}
\sum_{n=0}^{\infty}\P(L_{n,n}=n)y^n&=\sum_{n=0}^{\infty}\P(L_{n,n}=n,S_n< x_2)y^n+\sum_{n=0}^{\infty}\P(L_{n,n}=n,S_n\ge x_2)y^n.
\end{aligned}$$
However,
$$\sum_{n=0}^{\infty}\P(L_{n,n}=n,S_n< x_2)\leq x_2 \sum_{n=0}^{\infty}\P(L_{n,n}=n,S_n=0)<\infty,$$
so when $y\rightarrow1-$,
$$\sum_{n=0}^{\infty}\P(L_{n,n}=n,S_n\ge x_2)y^n\sim\frac{C_{\rho}(y)}{(1-y)^{\rho}}.$$
Using Tauberian's theorem, we get
$$\lim_{n\rightarrow\infty}\frac{\Gamma(\rho+1)}{n^{\rho}C_{\rho}(1-1/n)}\sum_{k=1}^{n}\P(L_{n,n}=n,S_n\ge x_2)=1,$$
then we get the proof.

\bigskip

 \acks
This research is partially sponsored by Natural Science Foundation of China (11671145),  Natural Science Foundation of Shanghai (16ZR1409700), and the program of China Scholarships Council (No.201806145024).

\end{document}